\documentclass{amsart}

\usepackage{amsfonts,amssymb}
\usepackage[latin1]{inputenc}

\newtheorem{thm}{Theorem}[section]

\newtheorem{lemma}[thm]{Lemma}
\newtheorem{prop}[thm]{Proposition}

\numberwithin{equation}{section}

\newcommand{\QED}{\hfill $\square$\vspace{2mm}}

\renewcommand{\epsilon}{\varepsilon}

\newcommand{\bbC}{\mathbb C}

\newcommand{\bbT}{\mathbb T}

\newcommand{\bbN}{\mathbb N}
\newcommand{\bbZ}{\mathbb Z}

\begin{document}

\title[Convergent Birkhoff vs. analytic integrability]
{Convergence versus integrability in Birkhoff normal form}%

\author{Nguyen Tien Zung}
\address{Laboratoire Emile Picard,
Département de Mathématiques, Université Toulouse III}
\email{tienzung@picard.ups-tlse.fr}


\date{Final version, 02/Dec/2003}%
\begin{abstract}
We show that any analytically integrable Hamiltonian system near an
equilibrium point admits a convergent Birkhoff normalization. The proof is
based on a new, geometric approach to the problem.
\end{abstract}
\maketitle
\section{Introduction}

Among the fundamental problems concerning analytic (real or complex) Hamiltonian
systems near an equilibrium point, one may mention the following two:

1) {\it Convergent Birkhoff}. In this paper, by ``convergent Birkhoff'' we mean the
existence of a convergent Birkhoff normalization, i.e. the existence of a local
analytic symplectic system of coordinates in which the Hamiltonian function will
Poisson commute with the semisimple part of its quadratic part.

2) {\it Analytic integrability}. By ``analytic integrability'' we mean the existence
of a complete set of local analytic functionally independent first integrals in
involution.

These problems have been studied by many classical and modern mathematicians,
including Poincar\'e, Birkhoff, Siegel, Moser, Bruno, etc. In this paper, we will be
concerned with the relations between the two problems. The starting point is that,
since both the Birkhoff normal form and the search for first integrals are a way to
simplify and solve Hamiltonian systems, these two problems must be very closely
related. Indeed, it has been known to Birkhoff \cite{Birkhoff} that, for nonresonant
Hamiltonian systems, convergent Birkhoff implies analytic integrability. The inverse
is also true, though much more difficult to prove \cite{Ito1}. What has been known
to date concerning ``convergent Birkhoff vs. analytic integrability'' may be
summarized in the following list. Denote by $q$ ($q \geq 0$) the {\it degree of
resonance} (see Section \ref{section:preliminary} for a definition) of an analytic
Hamiltonian system at an equilibrium point. Then we have :

a) When $q=0$ (i.e. for non-resonant systems), convergent Birkhoff is
equivalent to analytic integrability. The part ``convergent Birkhoff
implies analytic integrability'' is straightforward. The inverse has been
a difficult problem. Under an additional nondegeneracy condition involving
the momentum map, it was first proved by R\"ussmann \cite{Russmann} in
1964 for the case with two degrees of freedom, and then by Vey \cite{Vey}
in 1978 for any number of degrees of freedom. Finally Ito \cite{Ito1} in
1989 solved the problem without any additional condition on the momentum
map.

b) When $q=1$ (i.e. for systems with a simple resonance), then convergent Birkhoff
is still equivalent to analytic integrability. The part ``convergent Birkhoff
implies analytic integrability'' is again obvious. The inverse has been proved some
years ago by Ito \cite{Ito3} and Kappeler, Kodama and N\'emethi \cite{KKN}.

c) When $q \geq 2$ then convergent Birkhoff {\it does not} imply analytic
integrability. The reason is that the Birkhoff normal form in this case
will give us $(n-q+1)$ first integrals in involution, where $n$ is the
number of degrees of freedom, but additional first integrals don't exist
in general, not even formal ones. (A counterexample can be found in
Duistermaat \cite{Duistermaat112}, see also Verhulst \cite{Verhulst} and
references therein). The question ``does analytic integrability imply
convergent Birkhoff'' when $q \geq 2$ has remained open until now. The
powerful analytical techniques, which are based on the fast convergent
method and used in \cite{Ito1,Ito3,KKN}, could not have been made to work
in the case with non-simple resonances.

The main purpose of this paper is to complete the above list, by giving a positive
answer to the last question.

\begin{thm}
\label{thm:main} Any real (resp., complex) analytically integrable Hamiltonian
system in a neighborhood of an equilibrium point on a symplectic manifold admits a
real (resp., complex) convergent Birkhoff normalization at that point.
\end{thm}

An important consequence of Theorem \ref{thm:main} is that we may classify
degenerate singular points of analytic integrable Hamiltonian systems by
their analytic Birkhoff normal forms (see, e.g., \cite{ZungHDR} and
references therein).

The proof given in this paper of Theorem \ref{thm:main} works for any
analytically integrable system, regardless of its degree of resonance. Our
proof is based on a geometrical method involving homological cycles,
period integrals, and torus actions, and it is completely different from
the analytical one used in \cite{Ito1,Ito3,KKN}. In a sense, our approach
is close to that of Eliasson \cite{Eliasson-NF1990}, who used torus
actions to prove the existence of a smooth Birkhoff normal form for smooth
integrable systems with a nondegenerate elliptic singularity. The role of
torus actions is given by the following proposition (see Proposition
\ref{prop:torus2} for a more precise formulation):

\begin{prop}
\label{prop:torus} The existence of a convergent Birkhoff normalization is
equivalent to the existence of a local Hamiltonian torus action which preserves the
system.
\end{prop}

We also have the following result, which implies that it is enough to prove Theorem
\ref{thm:main} in the complex analytic case :

\begin{prop}
\label{prop:real}
A real analytic Hamiltonian system near an equilibrium point
admits a real convergent Birkhoff normalization if and only if it admits a complex
convergent Birkhoff normalization.
\end{prop}

Both Proposition \ref{prop:torus} and Proposition \ref{prop:real} are very
simple and natural. They are often used implicitly, but they have not been
written explicitly anywhere in the literature, to our knowledge.

The rest of this paper is organized as follows: In Section
\ref{section:preliminary} we introduce some necessary notions, and prove
the above two propositions. In Section \ref{section:proof} we show how to
find the required torus action in the case of integrable Hamiltonian
systems, by searching 1-cycles on the local level sets of the momentum
map, using an approximation method based on the existence of a formal
Birkhoff normalization and  \L ojasiewicz inequalities. This section
contains the proof of our main theorem, modulo a lemma about analytic
extensions. This lemma, which may be useful in other problems involving
the existence of first integrals of singular foliations (see
\cite{ZungHDR}), is proved in Section \ref{section:extension}, the last
section.

\section{Preliminaries}
\label{section:preliminary}

Let $H : U \to {\mathbb K}$, where $K = {\mathbb R}$ (resp., $K = {\mathbb
C}$) be a real (resp., complex) analytic function defined on an open
neighborhood $U$ of the origin in the symplectic space $({\mathbb K}^{2n},
\omega = \sum_{j=1}^n dx_j \wedge dy_j)$. When $H$ is real, we will also
consider it as a complex analytic function with real coefficients. Denote
by $X_H$ the symplectic vector field of $H$:
\begin{equation}
i_{X_H} \omega = - dH.
\end{equation}
Here the sign convention is taken so that $\{H,F\} = X_H(F)$ for any function $F$,
where
\begin{equation}
\{H,F\} = \sum_{j=1}^n \frac{dH}{dx_j}\frac{dF}{dy_j} -
\frac{dH}{dy_j}\frac{dF}{dx_j}
\end{equation}
denotes the standard Poisson bracket.

Assume that $0$ is an equilibrium of $H$, i.e. $dH(0)=0$. We may also put $H(0) =
0$. Denote by
\begin{equation}
H = H_2 + H_3 + H_4 + \hdots
\end{equation}
the Taylor expansion of $H$, where $H_k$ is a homogeneous polynomial of degree $k$
for each $k \geq 2$. The algebra of quadratic functions on $({\mathbb
K}^{2n},\omega)$, under the standard Poisson bracket, is naturally isomorphic to the
simple algebra $sp(2n, {\mathbb K})$ of infinitesimal linear symplectic
transformations in ${\mathbb K}^{2n}$. In particular,
\begin{equation}
H_2 = H_{ss} + H_{nil},
\end{equation}
where $H_{ss}$ (resp., $H_{nil}$) denotes the semi-simple (resp.,
nilpotent) part of $H_2$.

For each natural number $k \geq 3$, the Lie algebra of quadratic functions
on ${\mathbb K}^{2n}$ acts linearly on the space of homogeneous
polynomials of degree $k$ on ${\mathbb K}^{2n}$ via the Poisson bracket.
Under this action, $H_2$ corresponds to a linear operator $G \mapsto
\{H_2,G\}$, whose semisimple part is $G \mapsto \{H_{ss}, G\}$. In
particular, $H_k$ admits a decomposition
\begin{equation}
H_k = - \{H_2, L_k\} + H'_k \, ,
\end{equation}
where $L_k$ is some element in the space of homogeneous polynomials of
degree $k$, and $H'_k$ is in the kernel of the operator $G \mapsto
\{H_{ss}, G\}$, i.e. $\{H_{ss}, H'_k\} = 0$. Denote by $\psi_k$ the
time-one map of the flow of the Hamiltonian vector field $X_{L_k}$. Then
$(x',y') = \psi_k (x,y)$ (where $(x,y)$, or also $(x_j,y_j)$, is a
shorthand for $(x_1,y_1,\hdots,x_n,y_n)$) is a symplectic transformation
of $({\mathbb K }^{2n}, \omega)$ whose Taylor expansion is
\begin{equation}
\begin{array}{l}
x_j' = x_j (\psi(x,y)) = x_j - \partial L_k / \partial y_j + O(k) , \\
y_j' = y_j (\psi(x,y)) = y_j + \partial L_k / \partial x_j + O(k) ,
\end{array}
\end{equation}
where $O(k)$ denotes terms of order greater or equal to $k$. Under the new local
symplectic coordinates $(x'_j,y'_j)$, we have
$$
\begin{array}{ll}
H & = H_2(x,y) + \hdots + H_k(x,y) + O(k+1) \\
& = H_2(x'_j + \partial L_k / \partial y_j, y'_j - \partial L_k / \partial x_j) +
H_3(x'_j, y'_j) + \hdots + H_k(x'_j,y'_j) + O(k+1) \\
& = H_2(x'_j,y'_j) - X_{L_k}(H_2) + H_3(x'_j, y'_j) + \hdots + H_k(x'_j,y'_j) + O(k+1) \\
& = H_2(x'_j,y'_j) + H_3(x'_j, y'_j) + \hdots + H_{k-1}(x'_j,y'_j) + H'_k
(x'_j,y'_j) + O(k+1) .
\end{array}
$$

In other words, the local symplectic coordinate transformation $(x',y') =
\psi_k (x,y)$ of ${\mathbb K}^{2n}$ changes the term $H_k$ to the term
$H'_k$ satisfying $\{H_{ss},H'_k\}=0$ in the Taylor expansion of $H$, and
it leaves the terms of order smaller than $k$ unchanged. By induction, one
finds a series of local analytic symplectic transformations $\phi_k$
($k\geq 3$) of type
\begin{equation}
\phi_k(x,y) = (x,y) + \, terms \,\, of \,\, order \geq k-1
\end{equation}
such that for each $m \geq 3$, the composition
\begin{equation}
\Phi_m = \phi_m \circ \hdots \circ \phi_3
\end{equation}
is a symplectic coordinate transformation which changes all the terms of
order smaller or equal to $k$ in the Taylor expansion of $H$ to terms that
commute with $H_{ss}$.

By taking limit $m \to \infty$, we get the following classical result due
to Birkhoff et al. (see, e.g., \cite{Birkhoff,Bruno,SiMo}) :

\begin{thm}[Birkhoff et al.]
\label{thm:Birkhoff}
For any real (resp., complex) Hamiltonian system $H$ near an
equilibrium point with a local real (resp., complex) symplectic system of
coordinates $(x,y)$, there exists a formal real (resp., complex) symplectic
transformation $(x',y') = \Phi(x,y)$ such that in the coordinates $(x',y')$ we have
\begin{equation}
\label{equation:Birkhoff} \{H,H_{ss}\} = 0 ,
\end{equation}
where $H_{ss}$ denotes the semisimple part of the quadratic part of $H$.
\QED
\end{thm}

When Equation (\ref{equation:Birkhoff}) is satisfied, one says that the
Hamiltonian $H$ is in {\it Birkhoff normal form}, and the symplectic
transformation $\Phi$ in Theorem \ref{thm:Birkhoff} is called a {\it
Birkhoff normalization}. Birkhoff normal form is one of the basic tools in
Hamiltonian dynamics, and it has already been used in the 19th century by
Delaunay \cite{Delaunay} and Linstedt \cite{Lindstedt} for some problems
of celestial mechanics.

When a Hamiltonian function $H$ is in  normal form, then its first
integrals are also normalized simultaneously to some extent. More
precisely, one has the following folklore lemma, whose proof is
straightforward (see, e.g., \cite{Ito1,Ito3,KKN}) :

\begin{lemma}
\label{lemma:simultaneous} If $\{H_{ss},H\} = 0$, i.e. $H$ is in Birkhoff
normal form, and $\{H,F\} = 0$, i.e. $F$ is a first integral of $H$, then
we also have $\{H_{ss},F\} = 0 .$ \QED
\end{lemma}

Recall that the simple Lie algebra $sp(2n, {\mathbb C})$ has only one
Cartan subalgebra up to conjugacy. In terms of quadratic functions, there
is a complex linear canonical system of coordinates $(x_j,y_j)$ of
${\mathbb C}^{2n}$ in which $H_{ss}$ can be written as
\begin{equation}
\label{equation:Hs} H_{ss} = \sum_{i=1}^n \gamma_j x_jy_j
\end{equation}
where $\gamma_j$ are complex coefficients, called {\it frequencies}. (The
quadratic functions $\nu_1 = x_1y_1, \hdots,\nu_n = x_ny_n$ span a Cartan
subalgebra). The frequencies $\gamma_j$ are complex numbers uniquely
determined by $H_{ss}$ up to a sign and a  permutation. The reason why we
choose to write $x_jy_j$ instead of ${1 \over 2}(x_j^2+y_j^2)$ in Equation
(\ref{equation:Hs}) is that this way monomial functions will be
eigenvectors of $H_{ss}$ under the Poisson bracket:
\begin{equation}
\{ H_{ss}, \prod_{j=1}^n x_j^{a_j} y_j^{b_j} \} = (\sum_{j=1}^n (b_j -
a_j)\gamma_j) \prod_{j=1}^n x_j^{a_j} y_j^{b_j} .
\end{equation}
In particular, $\{H,H_{ss}\} = 0$ if and only if every monomial term
$\prod_{j=1}^n x_j^{a_j} y_j^{b_j}$ with a non-zero coefficient in the
Taylor expansion of $H$ satisfies the following relation, called a {\it
resonance relation}:
\begin{equation}
\sum_{j=1}^n (b_j - a_j)\gamma_j = 0 .
\end{equation}

In the nonresonant case, when there are no resonance relations except the
trivial ones, the Birkhoff normal condition $\{H,H_{ss}\} = 0$ means that
$H$ is a function of $n$ variables $\nu_1 = x_1y_1, \hdots, \nu_n =
x_ny_n$, implying complete integrability. Thus any nonresonant Hamiltonian
system is formally integrable \cite{Birkhoff,SiMo}.

More generally, denote by ${\mathcal R} \subset {\mathbb Z}^n$ the
sublattice of ${\mathbb Z}^n$ consisting of elements $(c_j) \in \bbZ^n$
such that $\sum c_j \gamma_j = 0$. The dimension of $\mathcal R$ over
${\mathbb Z}$, denoted by $q$, is called the {\it degree of resonance} of
the Hamiltonian $H$. Let $\mu^{(n-q+1)}, \hdots, \mu^{(n)}$ be a basis of
the resonance lattice $\mathcal R$. Let $\rho^{(1)},\hdots,\rho^{(n)}$ be
a basis of ${\mathbb Z }^n$ such that $\sum_{j=1}^n \rho^{(k)}_j
\mu^{(h)}_j = \delta_{kh}$ ($= 0$ if $k \neq h$ and $= 1$ if $k=h$), and
set
\begin{equation}
F^{(k)}(x,y) = \sum_{j=1}^n \rho^{(k)}_j x_j y_j
\end{equation}
for $1 \leq k \leq n$. Then we have $ H_{ss} = \sum_{k=1}^{n-q} \alpha_k
F^{(k)} $ with no resonance relation among $\alpha_1,\hdots,\alpha_{n-q}$.
The equation $\{H_{ss}, H\} = 0$ is now equivalent to
\begin{equation}
\{F_k , H \} = 0 \, \, \forall\ k = 1,\hdots,n-q .
\end{equation}

What is so good about the quadratic functions $F^{(k)}$ is that each
$iF^{(k)}$ (where $i = \sqrt{-1}$) is a periodic Hamiltonian function,
i.e. its holomorphic Hamiltonian vector field $X_{iF^{(k)}}$ is periodic
with a real positive period (which is $2\pi$ or a divisor of this number).
In other words, if we write $X_{iF^{(k)}} = X_k + i Y_k$ , where $X_k = J
Y_k$ is a real vector field called the real part of $X_{iF^{(k)}}$ (i.e.
$X_k$ is a vector field of ${\mathbb C}^{2n}$ considered as a real
manifold; $J$ denotes the operator of the complex structure of ${\mathbb
C}^{2n}$), then the flow of $X_k$ in ${\mathbb C}^{2n}$ is periodic. Of
course, if $F$ is a holomorphic function on a complex symplectic manifold,
then the real part of the holomorphic vector field $X_F$ is a real vector
field which preserves the complex symplectic form and the complex
structure.

Since the periodic Hamiltonian functions $iF^{(k)}$ commute pairwise (in
this paper, when we say ``periodic'', we always mean with a real positive
period), the real parts of their Hamiltonian vector fields generate a
Hamiltonian action of the real torus ${\mathbb T}^{n-q}$ on $({\mathbb
C}^{2n}, \omega)$. (One may extend it to a complex torus $({\mathbb
C}^{\ast})^{n-q}$-action, ${\mathbb C}^{\ast} = {\mathbb C} \backslash
\{0\}$, but we will only use the compact real part of this complex torus).
If $H$ is in (analytic) Birkhoff normal form, it will Poisson-commute with
$F^{(k)}$, and hence it will be preserved by this torus action.

Conversely, if there is a Hamiltonian torus action of ${\mathbb T}^{n-q}$
in $({\mathbb C}^{2n}, \omega)$ which preserves $H$, then the equivariant
Darboux theorem (which may be proved by an equivariant version of the
Moser path method, see, e.g., \cite{CDM}) implies that there is a local
holomorphic canonical transformation of coordinates under which the action
becomes linear (and is generated by $iF^{(1)},\hdots,iF^{(n-q)}$). Since
this action preserves $H$, it follows that $\{H, H_{ss} \} =0 $. Thus we
have proved the following

\begin{prop}
\label{prop:torus2} With the above notations, the following two conditions are
equivalent:

i) There exists a holomorphic Birkhoff canonical transformation of coordinates
$(x',y') = \Phi (x,y) $ for $H$ in a neighborhood of $0$ in ${\mathbb C}^{2n}$.

ii) There exists an analytic Hamiltonian torus action of ${\mathbb
T}^{n-q}$, in a neighborhood of $0$ in ${\mathbb C}^{2n}$, which preserves
$H$, and whose linear part is generated by the Hamiltonian vector fields
of the functions $iF^{(k)} = i\sum \rho^{(k)}_j x_j y_j$,
$k=1,\hdots,n-q$. \QED
\end{prop}

{\bf \emph{Proof of Proposition \ref{prop:real}}}. Let $H$ be a real
analytic Hamiltonian function which admits a local complex analytic
Birkhoff normalization, we will have to show that $H$ admits a local real
analytic Birkhoff normalization. Let $A : {\mathbb T}^{n-q} \times
({\mathbb C}^{2n}, 0) \to ({\mathbb C}^{2n}, 0)$ be a Hamiltonian torus
action which preserves $H$ and which has an appropriate linear part, as
provided by Proposition \ref{prop:torus}. To prove Proposition
\ref{prop:real}, it suffices to linearize this action by a local real
analytic symplectic transformation.

Let $F$ be a holomorphic periodic Hamiltonian function generating a
$\bbT^1$-subaction of $A$. Denote by $F^{\ast}$ the function $F^{\ast} (z)
= \overline{F(\bar{z})}$, where $z \mapsto \bar{z}$ is the complex
conjugation in ${\mathbb C}^{2n}$. Since $H$ is real and $\{H,F\} = 0$, we
also have $\{H,F^{\ast}\} = 0$. It follows that, if $H$ is in complex
Birkhoff normal form, we will have $\{H_{ss},F^{\ast}\}=0$, and hence
$F^{\ast}$ is preserved by the torus ${\mathbb T}^{n-q}$-action.
$F^{\ast}$ is a periodic Hamiltonian function by itself (because $F$ is),
and due to the fact that $H$ is real, the quadratic part of $F^{\ast}$ is
a real linear combination of the quadratic parts of periodic Hamiltonian
functions that generate the torus ${\mathbb T}^{n-q}$-action. It follows
that $F^{\ast}$ must in fact be also the generator of an ${\mathbb
T}^1$-subaction of the torus ${\mathbb T}^{n-q}$-action. (Otherwise, by
combining the action of $X_{F^{\ast}}$ with the ${\mathbb
T}^{n-q}$-action, we would have a torus action of higher dimension than
possible). The involution $F \mapsto F^{\ast}$ gives rise to an involution
$t \mapsto \bar{t}$ in ${\mathbb T}^{n-q}$. The torus action is reversible
with respect to this involution and to the complex conjugation:
\begin{equation}
\label{equation:involution}
\overline{A(t,z)} = A(\bar{t},\bar{z})
\end{equation}

The above equation implies that the local torus ${\mathbb T}^{n-q}$-action may be
linearized locally by a real transformation of variables. Indeed, one may use the
following averaging formula
\begin{equation}
z' = z'(z) = \int_{{\mathbb T}^{n-q}} A_1 (-t, A(t,z)) d\mu ,
\end{equation}
where $t \in {\mathbb T}^{n-q}$, $z \in {\mathbb C}^{2n}$, $A_1$ is the
linear part of $A$ (so $A_1$ is a linear torus action), and $d\mu$ is the
standard constant measure on ${\mathbb T}^{n-q}$. The action $A$ will be
linear with respect to $z'$ : $z' (A(t,z)) = A_1 (t,z'(z))$. Due to
Equation (\ref{equation:involution}), we have that $\overline{z'(z)} =
z'(\overline{z})$, which means that the transformation $z \mapsto z'$ is
real analytic.

After the above transformation $z \mapsto z'$, the torus action becomes
linear; the symplectic structure $\omega$ is no longer constant in
general, but one can use the equivariant Moser path method to make it back
to a constant form (see, e.g., \cite{CDM}). In order to do it, one writes
$\omega - \omega_0 = d\alpha$ and considers the flow of the time-dependent
vector field $X_t$ defined by $ i_{X_t} (t\omega + (1-t)\omega_0) =
\alpha$, where $\omega_0$ is the constant symplectic form which coincides
with $\omega$ at point $0$. One needs $\alpha$ to be ${\mathbb
T}^{n-q}$-invariant and real. The first property can be achieved, starting
from an arbitrary real analytic $\alpha$ such that $d\alpha = \omega -
\omega_0$, by averaging with respect to the torus action. The second
property then follows from Equation (\ref{equation:involution}).
Proposition \ref{prop:real} is proved. \QED

\section{Local torus actions for integrable systems}
\label{section:proof}

{\bf \emph{Proof of theorem \ref{thm:main}}}. According to Proposition
\ref{prop:real}, it is enough to prove Theorem \ref{thm:main} in the
complex analytic case. In this section, we will do it
 by finding local Hamiltonian $\bbT^1$-actions
which preserve the momentum map of an analytically completely integrable
system. The Hamiltonian function generating such an action will be a first
integral of the system, called an {\it action function} (as in
``action-angle coordinates''). If we find $(n-q)$ such $\bbT^1$-actions,
then they will automatically commute and give rise to a Hamiltonian
${\mathbb T}^{n-q}$-action.

To find an action function, we will use the following period integral
formula, known as {\it Mineur-Arnold formula}:
$$
P = \int_{\gamma} \beta \, ,
$$
where $P$ denotes an action function, $\beta$ denotes a primitive 1-form
(i.e. $\omega = d\beta$ is the symplectic form), and $\gamma$ denotes an
1-cycle (closed curve) lying on a level set of the momentum map.

To show the existence of such 1-cycles $\gamma$, we will use an approximation
method, based on the existence of a formal Birkhoff normalization.

Denote by ${\bf G} = (G_1 = H,G_2,\hdots,G_n): ({\mathbb C}^{2n},0) \to
({\mathbb C}^{n},0)$ the holomorphic momentum map germ of a given complex
analytic integrable Hamiltonian system. Let $\epsilon_0 > 0$ be a small
positive number such that ${\bf G }$ is defined in the ball $\{z =
(x_j,y_j) \in {\mathbb C}^{2n}, |z| < \epsilon_0\}$. We will restrict our
attention to what happens inside this ball. As in the previous section, we
may assume that in the symplectic coordinate system $z = (x_j,y_j) $ we
have
\begin{equation}
H = G_1 = H_{ss} + H_{nil} + H_3 + H_4 + \hdots
\end{equation}
with
\begin{equation}
H_{ss} = \sum_{k=1}^{n-q} \alpha_k F^{(k)}, \, \, F^{(k)} = \sum_{j=1}^n
\rho^{(k)}_jx_jy_j ,
\end{equation}
with no resonance relations among $\alpha_1,\hdots,\alpha_{n-q}$. We will
fix this coordinate system $z = (x_j,y_j)$, and all functions will be
written in this coordinate system.

The real and imaginary parts of the Hamiltonian vector fields of
$G_1,\hdots,G_n$ are in involution and define an {\it associated singular
foliation} in the ball $\{z = (x_j,y_j) \in {\mathbb C}^{2n}, |z| <
\epsilon_0\}$. Hereafter the norm in $\bbC^n$ is given by the standard
Hermitian metric with respect to the coordinate system $(x_j,y_j)$.
Similarly to the real case, the leaves of this foliation are called local
{\it orbits of the associated Poisson action}; they are complex isotropic
submanifolds, and generic leaves are Lagrangian and have complex dimension
$n$. For each $z$ we will denote the leaf which contains $z$ by $M_{z}$.
Recall that the momentum map is constant on the orbits of the associated
Poisson action. If $z$ is a point such that ${\bf G} (z)$ is a regular
value for the momentum map, then $M_z$ is a connected component of ${\bf
G}^{-1}({\bf G} (z))$.

Denote by
\begin{equation}
S = \{ z \in {\mathbb C}^{2n}, |z| < \epsilon_0, dG_1 \wedge dG_2 \wedge
\hdots \wedge dG_n (z) = 0 \}
\end{equation}
the singular locus of the momentum map, which is also the set of singular
points of the associated singular foliation. What we need to know about
$S$ is that it is analytic and of codimension at least 1, though for
generic integrable systems $S$ is in fact of codimension 2. In particular,
we have the following \L ojasiewicz inequality (see \cite{Lojasiewicz}) :
there exist a positive number $N$ and a positive constant $C$ such that
\begin{equation}
\label{equation:Lojasiewicz} |dG_1 \wedge \hdots \wedge dG_n (z) | > C
(d(z, S))^N
\end{equation}
for any $z$ with $|z| < \epsilon_0$, where the norm applied to $dG_1
\wedge \hdots \wedge dG_n (z)$ is some norm in the space of $n$-vectors,
and $d(z,S)$ is the distance from $z$ to $S$ with respect to the Euclidean
metric.

We will choose an infinite decreasing series of small numbers $\epsilon_m$
($m=1,2,\hdots$), as small as needed, with $\lim_{m \to \infty} \epsilon_m
= 0$, and define the following open subsets $U_m$ of ${\mathbb C}^{2n}$:
\begin{equation}
U_m = \{z \in {\mathbb C}^{2n}, |z| < \epsilon_m, d(z,S) > |z|^m  \}
\end{equation}

We will also choose two infinite increasing series of natural numbers
$a_m$ and $b_m$ ($m = 1,2,\hdots)$, as large as needed, with $\lim_{m \to
\infty} a_m = \lim_{m \to \infty} b_m = \infty$. It follows from
Birkhoff's Theorem \ref{thm:Birkhoff} and Lemma \ref{lemma:simultaneous}
that there is a series of local holomorphic symplectic coordinate
transformations $\Phi_m$, $m \in {\mathbb N}$, such that the following two
conditions are satisfied :

a) The differential of $\Phi_m$ at $0$ is identity for each $m$, and for any two
numbers $m,m'$ with $m' > m$ we have
\begin{equation}
\Phi_{m'}(z) = \Phi_m(z) + O (|z|^{a_m}) .
\end{equation}
In particular, there is a formal limit $\Phi_{\infty} = \lim_{m \to \infty} \Phi_m$.

b) The momentum map is normalized up to order $b_m$ by $\Phi_m$. More
precisely, the functions $G_j$ can be written as
\begin{equation}
G_j (z) = G_{(m)j} (z) + O (|z|^{b_m}), \, j=1,\hdots, n,
\end{equation}
with $G_{(m)j}$ such that
\begin{equation}
\{ G_{(m)j}, F_{(m)}^{(k)} \} = 0 \,\,\,\,  \forall j=1,\hdots, n, \,
k=1,\hdots,n-q .
\end{equation}
Here the functions $F_{(m)}^{(k)}$ are quadratic functions
\begin{equation}
F_{(m)}^{(k)} (x,y) = \sum_{j=1}^n \rho^{(k)}_j x_{(m)j} y_{(m)j}
\end{equation}
in local symplectic coordinates
\begin{equation}
(x_{(m)},y_{(m)}) = \Phi_m (x,y) .
\end{equation}

Notice that $F_{(m)}^{(k)}$ is a quadratic function in the coordinate
system $(x_{(m)},y_{(m)})$. But from now on we will use only the original
coordinate system $(x,y)$. Then $F_{(m)}^{(k)}$ is not a quadratic
function in $(x,y)$ in general, and the quadratic part of $F_{(m)}^{(k)}$
is $F^{(k)}$. The norm in $\bbC^{2n}$ used the estimations in this section
will be given by the standard Hermitian metric with respect to the
original coordinate system $(x,y)$.

Denote by $\gamma^{(k)}_m(z)$ the orbit of the real part of the periodic
Hamiltonian vector field $X_{iF_{(m)}^{(k)}}$ which goes through $z$. Then
for any $z' \in \gamma^{(k)}_m(z)$ we have $G_{(m)j}(z') = G_{(m)j}(z)$,
and $|z'| \simeq |z|$, i.e. $\lim_{z \to 0} {|z'| \over |z|} = 1$. (The
reason is that real part of the linear periodic Hamiltonian vector field
$X_{iF^{(k)}}$ also preserves the Hermitian metric of $\bbC^{2n}$, and the
linear part of $X_{iF_{(m)}^{(k)}}$ is $X_{iF^{(k)}}$). As a consequence,
we have
\begin{equation}
\label{equation:G} |{\bf G}(z') - {\bf G}(z)| = O (|z'|^{b_m}) .
\end{equation}
Note that, for each $m \in \bbN$, we can choose the numbers $a_m$ and
$b_m$ first, then choose the radius $\epsilon_m = \epsilon_m(a_m,b_m)$
sufficiently small so that the equivalence $O (|z'|^{b_m}) \simeq O
(|z|^{b_m})$ makes sense for $z \in U_m$.

On the other hand, we have
\begin{equation}
\begin{array}{l}
 |d G_1 (z') \wedge \hdots \wedge dG_n(z')| \\
 = |d G_{(m)1} (z') \wedge \hdots \wedge d G_{(m)n}(z')| + O (|z|^{b_m -1}) \\
 \simeq |d G_{(m)1} (z) \wedge \hdots \wedge d G_{(m)n}(z)| + O (|z|^{b_m -1}) \\
 = |d G_1 (z) \wedge \hdots \wedge dG_n(z)| + O (|z|^{b_m -1})
\end{array}
\end{equation}

We can assume that $b_m -1 > N$. Then for $|z| < \epsilon_m$ small enough,
the above inequality may be combined with \L ojasiewicz inequality
(\ref{equation:Lojasiewicz}) to yield
\begin{equation}
|d G_1 (z') \wedge \hdots \wedge dG_n(z')| > C_1 d(z,S)^N
\end{equation}
where $C_1 = C/2$ is a positive constant (which does not depend on $m$).

If $z \in U_m$, and assuming that $\epsilon_m$ is small enough, we have $d(z,S) >
|z|^m$, which may be combined with the last inequality to yield :
\begin{equation}
\label{equation:dG}
 |d G_1 (z') \wedge \hdots \wedge dG_n(z')| > C_1 |z|^{mN}
\end{equation}

Assuming that $b_m$ is much larger than $mN$, we can use the implicit function
theorem to project the curve $\gamma_m^{(k)}(z)$ on $M_z$ as follows :

For each point $z' \in \gamma_m^{(k)}(z)$, let $D_m(z')$ be the complex
$n$-dimensional disk centered at $z'$, which is orthogonal to the kernel
of the differential of the momentum map $G$ at $z'$, and which has radius
equal to $|z'|^{2mN}$. Since the second derivatives of $G$ are locally
bounded by a constant near $0$, it follows from the definition of
$D_m(z')$ that we have, for $|z| < \epsilon_m$ small enough :

\begin{equation}
\label{equation:dG2}
 |D{\bf G}(w) - D{\bf G}(z')| < |z|^{3mN/2} \,\, \forall w \in  D_m(z')
\end{equation}
where $D{\bf G}(w)$ denotes the differential of the momentum map at $w$,
considered as an element of the linear space of $2n \times n$ matrices.

Inequality (\ref{equation:dG}) together with Inequality
(\ref{equation:dG2}) imply that the momentum map ${\bf G}$, when
restricted to $D_m(z')$, is a diffeomorphism from $D(z')$ to its image,
and the image of $D_m(z')$ in ${\mathbb C}^{n}$ under ${\bf G}$ contains a
ball of radius $|z|^{4mN}$. (Because $4mN > 2mN + mN$, where $2mN$ is the
order of the radius of $D_m(z')$, and $mN$ is a majorant of the order of
the norm of the differential of $\bf G$. The differential of $\bf G$ is
``nearly constant'' on $D_m(z')$ due to Inequality (\ref{equation:dG2})).
Thus, if $b_m > 5mN$ for example, then Inequality (\ref{equation:G})
implies that there is a unique point $z''$ on $D_m(z')$ such that ${\bf
G}(z'') = {\bf G}(z)$. The map $z' \mapsto z''$ is continuous, and it maps
$\gamma^{(k)}_m(z)$ to some close curve $\tilde\gamma^{(k)}_m(z)$, which
must lie on $M_z$ because the point $z$ maps to itself under the
projection. When $b_m$ is large enough and $\epsilon_m$ is small enough,
then $\tilde\gamma^{(k)}_m(z)$ is a smooth curve with a natural
parametrization inherited  from the natural parametrization of
$\gamma^{(k)}_m(z)$, it has bounded derivative (we can say that its
velocity vectors are uniformly bounded by 1), and it depends smoothly on
$z \in U_m$.

Define the following action function $P_m^{(k)}$ on $U_m$ :

\begin{equation}
P_m^{(k)} (z) = \int_{\tilde\gamma_m^{(k)}(z)} \beta \ ,
\end{equation}
where $\beta = \sum x_j dy_j$ (so that $d\beta = \sum dx_j \wedge dy_j$ is the
standard symplectic form). This function has the following properties:

i) Because the 1-form $\beta = \sum x_j dy_j$ is closed on each leaf of
the Lagrangian foliation of the integrable system in $U_m$, $P_m^{(k)}$ is
a holomorphic first integral of the foliation. (This fact is well-known in
complex geometry: period integrals of holomorphic $k$-forms, which are
closed on the leaves of a given holomorphic foliation, over $p$-cycles of
the leaves, give rise to (local) holomorphic first integrals of the
foliation). The functions $P_m^{(1)}, \hdots, P_m^{(n-q)}$ Poisson commute
pairwise, because they commute with the momentum map.

ii) $P_m^{(k)}$ is uniformly bounded by 1 on $U_m$, because
$\tilde\gamma^{(k)}_m(z)$ is small together with its first derivative.

iii) Provided that the numbers $a_m$ are chosen large enough, for any $m' > m$ we
have that $P_m^{(k)}$ coincides with $P_{m'}^{(k)}$ in the intersection of $U_m$
with $U_{m'}$. To see this important point, recall that we have
\begin{equation}
P^{(k)}_m = P^{(k)}_{m'} + O(|z|^{a_m})
\end{equation}
by construction, which implies that the curve $\gamma^{(k)}_{m'}(z)$ is
$|z|^{a_m-2}$-close to the curve $\gamma^{(k)}_m(z)$ in $C^1$-norm. If $a_m$ is
large enough with respect to $mN$ (say $a_m > 5mN$), then it follows that the
complex $n$-dimensional cylinder
\begin{equation}
V_{m'}(z) = \{ w \in {\mathbb C}^{2n} \, | \, d(w,\gamma^{(k)}_{m'}(z)) < |z|^{2m'N}
\} \bigcap M_z
\end{equation}
lies inside (and near the center of) the complex $n$-dimensional cylinder
\begin{equation}
V_{m}(z) = \{ w \in {\mathbb C}^{2n} \, | \, d(w,\gamma^{(k)}_{m}(z)) < |z|^{2mN} \}
\bigcap M_z .
\end{equation}
On the other hand, one can check that $\tilde\gamma^{(k)}_m(z)$ is a retract of
$V_{m}(z)$ in $M_z$, and the same thing is true for the index $m'$. It follows
easily that $\tilde\gamma^{(k)}_{m'}(z)$ must be homotopic to
$\tilde\gamma^{(k)}_m(z)$ in $M_z$, implying that $P_m^{(k)}(z)$ coincides with
$P_{m'}^{(k)}(z)$.

iv) Since $P_m^{(k)}$ coincides with $P_{m'}^{(k)}$ in $U_m \bigcap U_{m'}$, we may
glue these functions together to obtain a holomorphic function, denoted by
$P^{(k)}$, on the union $U = \bigcup_{m=1}^{\infty} U_m$. Lemma
\ref{lemma:extension} in the following section shows that if we have a bounded
holomorphic function in $U = \cup_{m=1}^{\infty} U_m$ then it can be extended to a
holomorphic function in a neighborhood of $0$ in ${\mathbb C}^{2n}$. Thus our action
functions $P^{(k)}$ are holomorphic in a neighborhood of $0$ in ${\mathbb C}^{2n}$.

v) $P^{(k)}$ is a local periodic Hamiltonian function whose quadratic part is
$iF^{(k)} = i\sum \rho^{(k)}_j x_jy_j$. To see this, remark that
\begin{equation}
iF^{(k)}_m (z) = i\sum \rho^{(k)}_{j} x_{(m)j}y_{(m)j} =
\int_{\gamma_m^{(k)}(z)} \beta \ ,
\end{equation}
for $z \in U_m$. Since the curve $\tilde\gamma_m^{(k)}(z)$ is $|z|^{3mN}$-close to
the curve $\gamma_m^{(k)}(z)$ by construction (provided that $b_m > 4mN$), we have
that
\begin{equation}
P^{(k)} (z) = iF^{(k)}_m (z) + O(|z|^{3mN})
\end{equation}
for $z \in U_m$. Due to the nature of $U_m$ (almost every complex line in ${\mathbb
C}^{2n}$ which contains the origin $0$ intersects with $U_m$ in an open subset (of
the line) which surrounds the point $0$), it follows from the last estimation that
in fact the coefficients of all the monomial terms of order $< 3mN$ of $P^{(k)}$
coincide with that of $iF^{(k)}_m$, i.e. we have
\begin{equation}
P^{(k)} (z) = iF^{(k)}_m (z) + O(|z|^{3mN})
\end{equation}
in a neighborhood of $0$ in ${\mathbb C}^{2n}$. In particular, we have
\begin{equation}
P^{(k)} = \lim_{m \to \infty} iF^{(k)}_m \ ,
\end{equation}
where the limit on the right-and side of the above equation is understood as the
formal limit of Taylor series, and the left-hand side is also considered as a Taylor
series. This is enough to imply that $P^{(k)}$ has $i\sum \rho^{(k)}_j x_jy_j$ as
its quadratic part , and that $P^{(k)}$ is a periodic Hamiltonian of period $2\pi$
because each $iF^{(k)}_m$ is so. (If a local holomorphic Hamiltonian vector field
which vanishes at $0$ is formally periodic then it is periodic).

Now we can apply Proposition \ref{prop:torus2} and Proposition
\ref{prop:real} to finish the proof of Theorem \ref{thm:main}. \QED

\section{Holomorphic extension of action functions}
\label{section:extension}

The following lemma shows that the action functions $P^{(k)}$ constructed in the
previous section can be extended holomorphically in a neighborhood of $0$.

\begin{lemma}
\label{lemma:extension} Let $U = \bigcup_{m=1}^{\infty} U_m $, with $U_m = \{ x \in
{\mathbb C}^n, |x| < \epsilon_m, d(x,S) > |x|^m \}$, where $\epsilon_m$ is an
arbitrary series of positive numbers and $S$ is a local proper complex analytic
subset of ${\mathbb C}^n$  ($codim_{\mathbb C} S \geq 1$). Then any bounded
holomorphic function on $U$ has a holomorphic extension in a neighborhood of $0$ in
${\mathbb C}^n$.
\end{lemma}

{\bf \emph{Proof}}. Though we suspect that this lemma should have been
known to specialists in complex analysis, we could not find it in the
literature, so we will provide a proof here. When $n=1$ the lemma is
obvious, so we will assume that $n \geq 2$. Without loss of generality, we
can assume that $S$ is a singular hypersurface. We divide the lemma into
two steps :

{\it Step 1}. The case when $S$ is contained in the union of hyperplanes
$\bigcup_{j=1}^n \{x_j=0\}$ where $(x_1,\hdots,x_n)$ is a local
holomorphic system of coordinates. Clearly, $U$ contains a product of
non-empty annuli $\eta_j<|x_j|<\eta'_j$, hence $f$ is defined by a Laurent
series in $x_1,\cdots,x_n$ there. We will study the domain of convergence
of this Laurent series, using the well-known fact that the domain of
convergence of a Laurent series is logarithmically convex. More precisely,
denote by $\pi$ the map $(x_1,\cdots,x_n)\mapsto(\log |x_1|,\cdots,\log
|x_n|)$ from $({\mathbb C}^*)^n$ to ${\mathbb R}^n$, where ${\mathbb C}^*
= {\mathbb C} \backslash \{0\}$, and set
$$
E = \{ {\bf r} = (r_1,\hdots,r_n) \in {\mathbb R}^n \ \mid \ \pi^{-1}({\bf
r}) \subset U \}
$$
Denote by $Hull(E)$ the convex hull of $E$ in ${\mathbb R}^n$. Then since the
function $f$ is analytic and bounded in $\pi^{-1}(E)$, it can be extended to
abounded analytic function on $\pi^{-1}(Hull(E))$. On the other hand, by definition
of $U = \bigcup_{m=1}^{\infty} U_m$, there is a series of positive numbers $K_m$
(tending to infinity) such that $E \supset (\bigcup_{m=1}^{\infty} E_m)$, where
$$
E_m = \{(r_1,\hdots,r_n) \in {\mathbb R}^n \ \mid \ (r_j < -K_m\ \forall
j)\ ,\ (r_j > m r_i \ \forall j \ne i) \}
$$
It is clear that the convex hull of $\bigcup_{m=1}^{\infty} E_m$, with
each $E_m$ defined as above, contains a neighborhood of
$(-\infty,\hdots,-\infty)$, i.e. a set of the type
$$
\{(r_1,\hdots,r_n) \in {\mathbb R}^n \ \mid \ r_j < -K\ \forall j\} .
$$
It implies that the function $f$ can be extended to a bounded analytic function in
${\mathcal U} \bigcap ({\mathbb C}^*)^n$, where $\mathcal U$ is a neighborhood of
$0$ in ${\mathbb C}^n$. Since $f$ is bounded in ${\mathcal U} \bigcap ({\mathbb
C}^*)^n$, it can be extended analytically on the whole ${\mathcal U}$. Step 1 is
finished.

{\it Step 2.} Consider now the case with an arbitrary $S$. Then we can use
Hironaka's desingularization theorem \cite{Hironaka} to make it smooth.
The general desingularization theorem is a very hard theorem, but in the
case of a singular complex hypersurface a relatively simple constructive
proof of it can be found in \cite{BiMi-Resolution1999}. In fact, since the
exceptional divisor will also have to be taken into account, after the
desingularization process we will have a variety which may have normal
crossings. More precisely, we have the following commutative diagram
\begin{equation}
\begin{array}{ccccc}
Q & \subset & S' & \subset & M^n \\
\downarrow & & \downarrow & & \downarrow p \\
0 & \in & S & \subset & ({\mathbb C}^n,0)
\end{array} \ ,
\end{equation}
where $({\mathbb C}^n,0)$ denotes the germ of ${\mathbb C}^n$ at $0$ presented by a
ball which is small enough; $M^n$ is a complex manifold; the projection $p$ is
surjective, and injective outside the exceptional divisor; $S'$ denotes the union of
the exceptional divisor with the smooth proper submanifold of $M^n$ which is
desingularization of $S$ -- the only singularities in $S'$ are normal crossings; $Q
= p^{-1}(0)$ is compact. $M^n$ is obtained from $({\mathbb C}^n, 0)$ by a finite
number of blowing-ups along submanifolds.

Denote by $U' = p^{-1}(U)$ the preimage of $U$ under the projection $p$.
One can pull back $f$ from $U$ to $U'$ to get a bounded holomorphic
function on $U'$, denoted by $f'$. An important observation is that the
type of $U$ persists under blowing-ups along submanifolds.(Or
equivalently, the type of its complement, which may be called a
\emph{sharp-horn-neighborhood} of $S$ because it is similar to horn-type
neighborhoods of $S \setminus \{0\}$ used by singularists but it is sharp
of arbitrary order, is persistent under blowing-ups). More precisely, for
each point $x \in Q$, the complement of $U'$ in a small neighborhood of
$x$ is a ``sharp-horn-neighborhood'' of $S'$ at $x$. Since $S'$ only has
normal crossings, the pair $(U',S')$ satisfies the conditions of Step 1,
and therefore we can extend $f'$ holomorphically in a neighborhood of $x$
in $M^n$. Since $Q = p^{-1}(0)$ is compact, we can extend $f'$
holomorphically in a neighborhood of $Q$ in $M'$. One can now project this
extension of $f'$ back to $({\mathbb C}^n,0)$ to get a holomorphic
extension of $f$ in a neighborhood of $0$. The lemma is proved. \QED

{\it Remark}. The ``sharp-horn'' type of the complement of $U$ in the
above lemma is essential. If we replace $U$ by $U_m$ (for any given number
$m$) then the lemma is false.
\\

{\bf Acknowledgements}. I would like to thank Jean-Paul Dufour for
proofreading this paper, Jean-Claude Sikorav for supplying me with the
above proof of Step 1 of Lemma \ref{lemma:extension}, and Alexandre Bruno
for some critical remarks. I'm also thankful to the referee for his
pertinent remarks which helped improve the presentation of this paper.

\bibliographystyle{amsplain}

\end{document}